\documentclass{article}
\usepackage{amsmath}
\usepackage{amsthm}
\usepackage{amssymb}
\newtheorem*{weil}{Theorem A}

\newtheorem{thm}{Theorem}[subsection]
\newtheorem{lem}[thm]{Lemma}
\newtheorem{cor}[thm]{Corollary}
\newtheorem{prop}[thm]{Proposition}
\newtheorem{claim}[thm]{Claim}

\numberwithin{equation}{subsection}

\newcommand{\wpet}{Weil-Petersson}
\def\tec{Teichm\"uller\ }
\newcommand{\grad}{\operatorname{grad}}

\newcommand{\g}{\gamma}

\newcommand{\tg}{T_{g,n}}
\newcommand{\otg}{\overline T_{g,n}}
\newcommand{\SM}{\mathcal M}
\newcommand{\SO}{\mathcal O}

\newcommand{\ov}{\overline}
\newcommand{\vp}{\varphi}

\def\wh{\widehat}
\def\<{\left<}                          
\def\>{\right>}

\begin{document}

\title{The \wpet\ Isometry Group}
\author{Howard Masur\footnote{Partially supported by 
NSF grant DMS 9803497}\and Michael
Wolf\footnote{Partially
supported by NSF grant DMS 9971563}}
\date{}
\maketitle

\section{Introduction}

 Let $F=F_{g,n}$ be a surface of genus $g$ with $n$ punctures. 
We assume $3g-3+n>1$ and that $(g,n) \ne (1,2)$.
The purpose of this paper is to prove, for the Weil-Petersson
metric on Teichmuller space $T_{g,n}$, the analogue of Royden's
famous result \cite{Ro} that every complex analytic isometry of
$T_{g,0}$ with  respect to the Teichmuller metric is induced by
an
element of the mapping class group.  
His proof involved
a study of the local geometry of the cotangent bundle to
Teichmuller space.   Royden's result
was
extended to general $T_{g,n}$ by Earle-Kra \cite {EK}, without
any smoothness assumption on the isometry and with  a stronger
 local result.  They showed that if $2g+n>4$ and $2g'+n'>4$, and
if $f$ is an isometry from an open set  $U\subset T_{g,n}$ to
$T_{g',n'}$, then $T_{g,n}
=T_{g',n'}$ and $f$ is the restriction of 
an isometry induced by
an element of the
extended mapping class group. 
Later Ivanov
\cite{Iv97} gave an alternative proof of Royden's
theorem 
based upon the asymptotic geometry of Teichmuller space and the
result  that  the group of automorphisms of the curve complex 
$C(F)$ (see below)
 coincides with the mapping class group. The automorphism result
was 
later extended to the cases of punctured surfaces of genus
$g\le1$ (with $(g,n) \ne (1,2)$)
by Korkmaz
\cite {Kor97}, and at the same time proved for
general $(g,n) \ne (1,2)$ by Luo \cite{Luo00}.

We prove
\begin{weil}\label{weil thm}
For $3g-3+n>1$ and $(g,n)\neq (1,2)$, every \wpet\ isometry of
\tec space $\tg$
is induced by an
element of the extended mapping class group $Mod^*(g,n)$.
\end{weil}

Our proof of this result is modelled somewhat on Ivanov's proof. 
We outline the ideas here.  It is
well-known that the Weil-Petersson metric is not complete
(\cite{Wol75} and \cite{Chu75}).
To complete the metric  one adds a frontier 
$A_{g,n}$ to $T_{g,n}$;  this
frontier consists of a union of lower dimensional Teichmuller
spaces.  Each such space  consists of Riemann surfaces  
with
nodes or punctures.  These surfaces are obtained by  pinching
nontrivial curves of $F$. 
Each  Teichmuller space on the frontier carries its own
 Weil-Petersson metric and with this Weil-Petersson metric, this 
Teichmuller space on the frontier is isometrically embedded in
the completion.  We extend the 
isometry to the completion and show that each Teichmuller space
on the frontier is preserved by the isometry. This 
self-map of the frontier then
induces an automorphism of the
complex of curves, $C(F)$.  By Ivanov's result 
(and the extensions \cite {Kor97} \cite{Luo00}
by Korkmaz and Luo to low genus) 
the automorphism
is induced by an element of the extended mapping
class group.  (Ivanov's theorem is  known not to hold
for $(g,n)=(1,2)$: see Luo \cite {Luo00}.  Luo showed that
$C(F_{1,2})$ 
is isomorphic to $C(F_{0,5})$; the automorphism group of this
curve
complex 
is the extended mapping class group $Mod^*(0,5)$ and yet
$Mod^*(1,2)$ 
is a subgroup of index $5$ in $Mod^*(0,5)$).

Thus our isometry induces an extended  mapping class group
element
at
``infinity''.  We then
show that the given isometry and the isometry induced by the
corresponding element of the extended mapping class group act
identically
on Teichmuller space.  The main tool, as in
Ivanov, is to study geodesics in the space.
In particular, after some formalities,
we find that it is enough to assume that the
isometry acts as the identity on the 
frontier $A_{g,n}$, and
we study the (totally geodesic submanifold)
$Fix$ which is fixed by the isometry.
General facts about fixed-point (proper)
subsets of isometries in CAT(0) spaces 
yield points in $T_{g,n}$ at arbitrary distance 
from $Fix$.  Yet, since \tec space is a space
of Riemann surfaces, additional estimates come from
considering the functions $l_c(x): T_{g,n} \to \Bbb{R}$,
defined as the hyperbolic length of the curve $c$ 
on the hyperbolic surface $x \in T_{g,n}$. A very strong
such estimate is due to Wolpert \cite{Wol87}, who shows
that such functions are convex along Weil-Petersson
geodesics. This allows us, following Wolpert, to give
a center-of-mass argument that shows that $Fix$ 
is non-empty.  Combining this with the proof of Wolpert's result
\cite{Wol76}
that the Weil-Petersson metric is not complete, we find that
$Fix$, along with all points in $T_{g,n}$, is uniformly close 
to a special set of frontier points.  This contradiction
with the general fact above about fixed-point proper subsets in
CAT(0) 
spaces shows that $Fix$ is all of $T_{g,n}$, proving the theorem.

We would like to thank the referees for many helpful suggestions
and for pointing out a gap in the original proof. 
We would like to thank Cliff Earle and Sumio Yamada 
for bringing this problem to our attention and to 
Cliff Earle, Jeff Brock,  and Amie Wilkinson for  helpful
suggestions. We also
appreciate some helpful remarks by Feng Luo about his 
paper.

\subsection{Teichmuller space, mapping class group  and the
Weil-Petersson metric}
Denote by $\SM$ the set of all smooth Riemannian metrics on $F$.
Choose an orientation for $F$, and define the set of all
similarly oriented complete hyperbolic
structures on $F$ by $\SM_{-1}$: here
$\SM_{-1}$ naturally includes in $\SM$.  
By the uniformization theorem,
$\SM_{-1}$ can be identified with the
set of all conformal structures on $F$, with the given
orientation.   Equivalently, this is the same as the set of all
complex structures or Riemann surface structures on $F$ with the
given orientation. The group
of orientation preserving diffeomorphisms   $Diff^+(F)$
acts on $\SM_{-1}$ by pull-back.  Let $Diff_0(F)$
the subgroup of
diffeomorphisms
isotopic to the identity.   The Teichmuller space 
$T_{g,n}$ is defined to be  $$T_{g,n}=\SM_{-1}/Diff_0(F).$$
We can equivalently define $T_{g,n}$ by fixing a complex
structure $S_0$ on $F$ and defining $T_{g,n}$ as the set of 
equivalence classes of pairs $(S,f)$ where $f:S_0\to S$ is a 
sense-preserving quasiconformal map from $S_0$ to $S$.  
Two pairs $(S,f)$ and $(S',f')$ are equivalent if there is a
conformal map $h:S\to S'$ such that $h\circ f$ is homotopic to
$f'$.  

The mapping class group $Mod(g,n)$ is defined to be
\begin{displaymath}Mod(g,n)=Diff^+(F)/Diff_0(F)
\end{displaymath}
 and the 
  extended mapping class group is defined  by
\begin{displaymath}Mod^*(g,n)=Diff(F))/Diff_0(F).
\end{displaymath}
The group $Mod^*(g,n)$ acts on $T_{g,n}$ as follows.   We may
choose $S_0$ so that it admits an antiholomorphic reflection
$j:S_0\to S_0$.  
Let $\Psi \in Mod(g,n)$ be
represented by $\psi:S_0\to S_0$.   
For $\psi:S_0\to S_0$ orientation preserving,
$\psi\cdot (S,f)=(S,f\circ \psi^{-1})$.  Any orientation
reversing
diffeomorphism
of $S_0$ can be expressed as $\psi\circ j$ for some orientation
preserving $\psi$.  Then  $(\psi\circ j)\cdot (S,f)$ is the point
$(S^*,f\circ j\circ \psi^{-1})$, where $S^*$ is the conjugate
Riemann
surface to $S$; that is, the coordinate charts of $S^*$ are those
of $S$ followed
by complex conjugation.  

The Moduli space $\SM_{g,n}$ is defined to be 
$\SM_{g,n}=T_{g,n}/Mod_{g,n}$.

It is well-known (see \cite{Nag})
that $T_{g,n}$ has a complex structure.  The cotangent space at a
point $X \in T_{g,n}$ 
is the space of holomorphic quadratic differentials $\phi(z)dz^2$
on the Riemann surface $X$. On $X$ there is a
pairing of quadratic differentials and  
 Beltrami differentials $\mu(z)\frac{d\ov z}{dz}$ 
(i.e. tensors of type $(-1,1)$) on $X$.   The
pairing is given  by $$<\mu,\phi>=Re \int_X
\mu(z)\phi(z)
dz\wedge d\overline z.$$
Infinitesimally trivial Beltrami differentials $\mu$ are ones
such that
$$<\mu,\phi>=0$$ for all holomorphic
$\phi$.  The tangent space at $X$ is the space of Beltrami
differentials modulo the infinitesimally trivial ones (see
\cite{Ahl}).

  The
Weil-Petersson co-metric on $T_{g,n}$ is defined by the $L\sp
2$-product on the cotangent
bundle 
$$<\phi\sb 1(z)dz\sp 2,\phi\sb 2(z)dz\sp 2>=\frac{\sqrt
{-1}}{2}\int\sb
\Sigma \frac{\phi\sb 1(z)\ov\phi\sb 2(z)}{\lambda(z)}dz\wedge
d\overline
z,$$ where $\lambda (z)\vert dz\vert \sp 2$
is the hyperbolic metric on $F$.
  The metric is then defined on the tangent space by duality;
alternatively, the metric is induced from the natural inner
product
on the tangent bundle $T\SM$ to $\SM$ along $\SM_{-1}$,
after projecting to $T_{g,n}$ (see \cite{Tro92}). 
 The metric is Kahler \cite{Ahl}.
The  major properties that we will use are that the 
metric has negative sectional curvature (\cite{Roy74}, 
\cite{Tro86}, \cite{Wol86}) and 
Wolpert's remarkable result \cite{Wol87}
that even though the metric is not complete, (see below) it is
geodesically
 convex: there is a geodesic joining any two points,  unique
because of the
 negative curvature. We also note that the 
action of $Mod^*(g,n)$ on $\tg$ is isometric  with
respect 
to the Weil-Petersson metric.

\subsection{The complex of curves, Ivanov's theorem, and the
frontier of Teichmuller space}

We define a complex $C(F)$ as follows.  The vertices of $C(F)$
are homotopy classes of
homotopically nontrivial, nonperipheral simple closed curves on
$F$.
An edge of $C(F)$  consists of a pair of homotopy classes of
disjoint simple closed
curves. More generally, a
$k$-simplex
consists
of $k+1$ homotopy classes of 
mutually disjoint simple closed curves.  The maximal
number of mutually disjoint
simple closed curves is $3g-3+n$ so that  $C(F)$ is a $3g-4+n$
dimensional simplicial  complex.
The extended mapping class group $Mod^*(g,n)$ clearly acts 
on $C(F)$ by simplicial automorphisms.  Namely for
$\psi\in
Mod^*(g,n)$ and
$v_k=\{\beta_1,\ldots,\beta_{k+1}\}$ a $k$-simplex, 
the image $\psi(v_k)$ is the $k$-simplex
 $\{\psi(\beta_1),\ldots,\psi(\beta_{k+1})\}$.  Ivanov proved
\cite{Iv97}
that in all but a few low genus cases, every simplicial 
automorphism of
$C(F)$
is induced by some $g\in Mod^*(g,n)$
in the extended mapping class group.  This was later extended by
Korkmaz \cite {Kor97} to all  cases  except $(g,n) = (1,2)$.
At the same time, Luo \cite{Luo00} gave an independent proof
of all cases $(g,n) \ne (1,2)$
(with the explanation of the case $(g,n) = (1,2)$).

\vskip3mm

It is well-known that $\tg$  is not compact.  
For a simple closed curve $\beta$ we can define a
function 
$l_\beta:\tg\to R$ by setting $l_\beta(x)$ to be the 
length, in the hyperbolic metric on $x$, 
of the geodesic in the homotopy class of $\beta$. For a
collection of curves $C=\beta_1,\ldots,\beta_N$, let 
\begin{displaymath}
l_C=\sum_{i=1}^N l_{\beta_i}
\end{displaymath}

One way to leave
all
compacta in $\tg$ is to choose a simplex $v_k$ in $C(F)$,
that is
a set of disjoint simple closed curves
$\beta_1,\ldots,\beta_{k+1}$, and form a sequence of Riemann
surfaces  along which the $l_{\beta_i}$ go to $0$. 
This motivates the definition of the  
{\it augmented Teichmuller space} (\cite {Be},\cite{Ab71}).

Specifically, associate to $v_k$ the (possibly disconnected)
surface
$S \setminus \{\beta_1,\ldots,\beta_{k+1}\}$ whose components are
punctured
surfaces $S_1,\ldots,S_p$.  
Each $S_i$ has its own Teichmuller space $T(S_i)$ and we let 
$\SO(v_k)$ be the product $T(S_1)\times\ldots\times 
T(S_p)$ of Teichmuller spaces $T(S_i)$.

An 
alternative description is given in terms of surfaces
with nodes (i.e. complex spaces in which each point has
a neighborhood homeomorphic to either $\{|z|< \epsilon\}$ 
(regular points)
or $\{(z,w) \in \Bbb{C}^2| zw=0, |z|< \epsilon, |w|< \epsilon\}$
(the {\it nodes})).  We identify each of those
curves $\beta_i$ to a point;  the resulting space, 
say $S(v_k)$, is 
homeomorphic to a surface with nodes and we 
set $\SO(v_k)= T(S(v_k))$. Here $T(S(v_k))$ is
defined to be the product \tec space of the punctured
surfaces obtained by removing the nodal points from
$S(v_k)$.

Thus, for each simplex $v_k$ in $C(F)$
we consider a  frontier Teichmuller space $\SO(v_k)$ 
which is of complex
dimension $3g-4+n-k$.  
  We denote by $A_{g,n}$ the union of all
these frontier Teichmuller spaces at infinity
and by  $\ov T_{g,n}=T_{g,n}\cup A_{g,n}$.    
There is a standard way (see \cite {Be},\cite{Ab71})
to topologize this union,
sometimes referred to as the augmented \tec space.

In this topology    $\ov T_{g,n}$ is
not compact, nor is it locally compact near a frontier point in
$\SO(v_k)$: 
take $\beta\in v_k$, and let $\tau$ be the  Dehn twist about
$\beta$.
Then the
 $\tau$ orbit of a
point $x\in T_{g,n}$  does not have a convergent subsequence. 

By way of contrast, the action of the mapping class
group
extends 
to an action on $\ov T_{g,n}$ and the quotient is a
compactification 
$\ov \SM_{g,n}$ of the moduli
space $\SM_{g,n}$, commonly called the 
{\it Deligne-Mumford compactification} (see \cite{Ab77},
\cite{Nag}).

We
may
think of $\ov T_{g,n}$ as a stratified space,  because if $v_l$
is a subsimplex of $v_k$, then the
Teichmuller
space $\SO(v_k)$  is on the frontier of the 
Teichmuller space $\SO(v_l)$.  This is because $k-l$
curves of the surfaces in 
$\SO(v_l)$ have lengths that have become $0$, or equivalently,
the surfaces have acquired an additional  $k-l$ nodes.  All of
these
spaces lie in $A_{g,n}$,  the union of the frontier spaces of
$T_{g,n}$. For each $\SO(v_k)$ we will denote by $Fr(\SO(v_k))$
the union of its frontier spaces. Now the frontier $A_{g,n}$
is connected, although if we fix  $k$, then each $\SO(v_k)$ is a
component of 
the union over all frontier Teichmuller spaces  of
that dimension.  

The case of a maximal simplex 
$v_{3g-4+n}$ is especially important, for in that case the
resulting 
 surface with punctures  is a union of thrice-punctured 
spheres.  Since the conformal structure (or equivalently,
hyperbolic structure)
on a three times
punctured sphere is unique, the corresponding frontier
Teichmuller space $\SO(v_{3g-4+n})$ 
are singletons.    
We call these maximally pinched frontier spaces.  They will play
a crucial role in the sequel.

\subsection{Incompleteness of Weil-Petersson metric and extension
of isometries}
Wolpert \cite{Wol75} and Chu \cite{Chu75} proved that the
Weil-Petersson metric is not complete
on 
$T_{g,n}$.   In fact they showed that the 
Weil-Petersson distance in Teichmuller space to any
frontier space $\SO(v_k)$ is finite. 
Thus  we can complete the metric by adding $A_{g,n}$,
inducing a metric $dist_{v_k}$ on each frontier 
Teichmuller space $\SO(v_k)$. Of course, each 
frontier Teichmuller space $\SO(v_k)$ 
already has its own  Weil-Petersson
metric, written $d_{v_k}(\cdot,\cdot)$.  
\begin{lem}
For two
points $p_0,p_1$ in the same frontier space $\SO(v_k)$, we have
$d_{v_k}(p_0,p_1) = dist_{v_k}(p_0,p_1)$.
\end{lem}
\begin{proof}
The Weil-Petersson metric tensor in $T_{g,n}$   extends
continuously 
to the Weil-Petersson metric tensor in $\SO(v_k)$.  (\cite{Ma76})
This implies that $d_{v_k}(p_0,p_1) \ge dist_{v_k}(p_0,p_1)$. 
On the  other hand, suppose a length-minimizing path 
in the completion metric joining
$p_0$ and $p_1$ enters $T_{g,n}$, and is thus a Weil-Petersson
geodesic $\sigma$ there. 
Without loss of generality we can assume $\sigma$  lies inside
$T_{g,n}$ except for its endpoints. Then  $l_{v_k}$ tends to 
zero near its endpoints, but is positive somewhere
in its interior. 
This
contradicts Wolpert's convexity result \cite{Wol87} which says
that the functions $l_\beta$ are strictly convex along Weil-
Petersson geodesics.  Thus  $d_{v_k}(p_0,p_1) \le
dist_{v_k}(p_0,p_1)$,
completing the proof.
\end{proof}

We will refer to $d_{WP}(\cdot,\cdot)$ as the completed metric on
$\ov T_{g,n}$.  
We emphasize that the restriction 
of $d_{WP}$ to any space $\SO(v_k)$ is the
Weil-Petersson metric on $\SO(v_k)$.

\vskip3mm


\noindent {\bf Remarks.} (i) Slightly stronger 
conclusions may also be drawn about this
situation with additional use of 
Wolpert's convexity result (\cite{Wol87}).  In particular,
we see that each frontier space $\SO(v_k)$, with the metric
$d_{v_k}$,
is geodesically embedded in $\otg$, in the sense that any
geodesic
connecting a pair of points in $\SO(v_k)$ lies entirely
in that component $\SO(v_k)$.  To see this fact (to our
knowledge,
first written down in \cite{Y}), note first that such a geodesic,
say $\sigma$, cannot meet a higher dimensional component
$\SO(v_{k-j})$,
for $v_{k-j} \subset v_{k}$; 
as above, this is because the length of any curve
in $v_{k}\setminus v_{k-j}$ is a non-negative convex function
which vanishes at its endpoints, hence vanishes identically.
On the other hand, since $\SO(v_k)$ is a product of 
Teichmuller spaces (with the Weil-Petersson metric)
of punctured surfaces, we see that on the geodesic $\sigma$,
the length functions are bounded by the maximum of their 
values at their endpoints, and hence $\sigma$ meets no
frontier spaces $\SO(v_{k+l})$ on the frontier of $\SO(v_k)$.

(ii) We note that since the Weil-Petersson metric has negative
curvature
and is geodesic convex, 
$\tg$ is a $CAT(0)$ space with this metric. It is a general fact 
(\cite{BH}, Corollary II.3.11)
that  the metric
completion of a geodesically convex
$CAT(0)$ space is $CAT(0)$.  Thus in particular $\otg$ is 
$CAT(0)$
and thus between any two points there is a unique geodesic. 
We however will not need to use this last fact. 

We will need the following lemmas. 

\begin{lem} 
Fix a frontier space $\SO(v_k)$ (which may be $\tg$ itself). Then,
for all $\epsilon$ sufficiently small, there exists
$L=L(\epsilon)$, depending only on $\epsilon$ 
with the following property. If 
$\SO(v_{3g-4+n})$ is any maximally pinched  frontier point on the
frontier of 
$\SO(v_k)$,   and $C$ is the
collection
of curves in $v_{3g-4+n}\setminus v_k$, then if
$d_{WP}(x,\SO(v_{3g-4+n}))\leq \epsilon$, we have $l_C(x)\leq L$.
\end{lem}
\begin{proof}
If not, there exists a sequence $\{x_j\}\in\SO(v_k)$ 
and a sequence $\{\SO(v_{{3g-4+n},j})\}$ of maximally
pinched frontier points
on the frontier of $\SO(v_k)$
with 
$d_{WP}(x_j,\SO(v_{{3g-4+n},j}))\to 0$ and $l_{C_j}(x_j)\to\infty$;
here, $C_j$ of course refers to the curves 
$v_{{3g-4+n},j}\setminus v_k$.
Since there are but a finite number of homotopy classes
of maximally pinched surfaces, we can find subsequences,
again called $\{x_j\}$ and $\{\SO(v_{{3g-4+n},j})\}$ and
a sequence $\{f_j\} \subset Mod(g,n)$ so that
$f_j\SO(v_{{3g-4+n},j})=\SO(w_{3g-4+n})$, where 
$\SO(w_{3g-4+n})$ is some single maximally pinched  frontier point on the
frontier of $\SO(v_k)$.
Then let $C$ denote the curve system $C=w_{3g-4+n}\setminus v_k$,
so that $C=f_jC_j$, and set $y_j=f_jx_j$.
Then, since $f_j$ induces a Weil-Petersson isometry,
we have $d(y_j, \SO(w_{3g-4+n})) = d(f_jx_j, f_j\SO(v_{{3g-4+n},j}))
=d(x_j, \SO(v_{{3g-4+n},j})) \to 0$, while
$l_C(y_j)=l_{f_jC_j}(f_jx_j)=l_{C_j}(x_j)\to \infty$: this
last equality follows from $f_j$ amounting to but 
a consistent relabelling of curves and hyperbolic surfaces.
 But then the first limit implies that the sequence
$\{y_j\}$ converges
to $\SO(w_{3g-4+n})$.  However this implies   $l_C(y_j)\to 0$ and
we have a contradiction with the second limit statement.  
\end{proof}
     
\begin{lem}
Fix a simplex $v_k \subset C(F)$.
Given $\rho>0$ and $M>0$ there is a $\delta=\delta(\rho,M)$ with
the following property. 
Let $C$ be a collection of curves on the surface $S= F \setminus
v_k$.
If $x\in\SO(v_k)$ satisfies $l_C(x)\leq M$ and
$d_{WP}(x,Fr(\SO(v_k))\geq\rho$,  
then  any $y$ in the Weil-Petersson $\delta$
ball
about $x$ satisfies
$l_C(y)\leq 2M$
\end{lem}

Here $Fr(\SO(v_k)$ denotes the frontier of $\SO(v_k)$.

Note there must be some condition on the location of 
$x$ for the conclusion to
hold. Take a  frontier space $\SO(w)$ of $\SO(v_k)$ obtained 
by pinching along a curve $\beta$ not in 
$C$, but which intersects some curve in $C$.   In any $\delta$
neighborhood of a point in $\SO(w)$
there are points  $x,y\in\SO(v_k)$ which differ
 by arbitrarily large powers of the Dehn twist about $\beta$. 
But then the ratio of $l_C(x)$ to  $l_C(y)$ can be made 
 arbitrarily large.

\begin{proof}
If the lemma were not true, there would be  sequences
$\{x_j\},\{x_j'\}\in\SO(v_k)$ and a sequence $\{C_j\}$ of curve families
 such that $d_{WP}(x_j,x_j')\to 0$, $d_{WP}(x_j,Fr(\SO(v_k))\geq
\rho$, 
 and $l_{C_j}(x_j)\leq M$, while $l_{C_j}(x_j')> 2M$.
>From the first two conditions   we can find 
subsequences again denoted $\{x_j\},\{x_j'\}$  and a sequence
$\{f_j\} \subset Mod(g,n)$
such that $\{y_j=f_j(x_j)\}$ and $\{y_j'=f_j(x_j')\}$  converge to the
same point
 $y_0\in\SO(v_k)$.  Fix $K$ a compact subset of $\SO(v_k)$
containing 
$y_j,y_j'$ for $j$ large.    
  Now 
$\{C_j'=f_j(C_j)\}$ is a sequence  of curve families such that
$l_{C_j'}(y_j)  = l_{f_j(C_j)}(f_j(x_j)) = l_{C_j}(x_j) \leq M$, 
while $l_{C_j'}(y_j')>2M$.  Since $K$ is compact, 
there are only
finitely many curves $\beta$ such that $l_\beta(z)\leq M$ for
some $z\in K$.    Thus $\bigcup_jC_j'$ is actually a finite set of
curves,
and by passing to subsequences we can assume the sequence $\{C_j'\}$
is a
fixed set $C_0$.  But now the length functions $l_\beta(\cdot)$ 
are continuous on $\SO(v_k)$ and
we have contradicted the fact that 
both $\{y_j\}$ and $\{y_j'\}$ converge to $y_0$.    
\end{proof}

Now suppose
$I$ is an isometry of $T_{g,n}$ into itself in the Weil-Petersson
metric.
Since we do not a priori assume that the isometry $I$ is
surjective,
we first show
\begin{lem}
 $ I:\tg\to \tg$ is surjective.
\end{lem}
\begin{proof}  Since $I$ is an
open map, we need to
      show that $I$
      is proper.  Suppose on the contrary that there is a
sequence $\{x_j\}$ leaving every
      compact set in
      $T_{g,n}$ such that $\{I(x_j)\}$ lies in a compact set $K$ of
$\tg$.  We first show that no subsequence of $\{x_j\}$ can  project to
a precompact open set $M$ in the moduli space $\SM_{g,n}$. 
 For as  the
compact 
set $K$ has finite Weil-Petersson diameter, so does 
 $I^{-1}(K)\subset T_{g,n}$; since a set of finite diameter
can intersect 
only finitely many disjoint balls of fixed diameter,  $I^{-1}(K)$ 
could intersect  but a finite number of preimages of the 
precompact open set $M\subset \SM_{g,n}$ under the projection map 
$T_{g,n}\to \SM_{g,n}$.  Thus the subsequence
$\{x_j = I^{-1}I(x_j)\}$ would be contained in the closure
of the union of those finite number of precompact
preimages, and hence
would be contained in a compact set in $T_{g,n}$, contrary to
hypothesis. 

We conclude that the entire sequence projects to a sequence that
leaves every compact subset of $\SM_{g,n}$.    
However along such  a sequence $\{x_j\}$
the infimum of  the scalar
curvature of the Weil-Petersson metric
is  $-\infty$
\cite{Tra92}, 
while
bounded      
in the compact set $K$.  But this is a contradiction to the fact
that an
      isometry preserves
      scalar
      curvature.   Thus $I$ is in fact surjective. 
\end{proof}

It is immediate that $I$ is therefore also invertible.

\vskip3mm

We discuss next how $I$ extends to an isometry
$\ov I:\otg\to\otg$ and that this extension preserves 
each frontier space $\SO(v_k)$.

\begin{lem}
 The isometry $I$ extends to an isometry 
$\ov I:\otg\to\otg$
of $\ov T_{g,n}$ which
is surjective. 
For each $k=0,\ldots,3g-4+n$, the isometry 
$\ov I$ sends each frontier space
$\SO(v_k)$ determined by a $k$-simplex $v_k$ onto a frontier
space
$\SO(w_k)$, where $w_k$ is another $k$-simplex. 
\end{lem}

We have the immediate corollary.
\begin{cor} $I$ induces an automorphism $\wh I$  of the curve
complex
$C(F)$.
\end{cor}
We prove Lemma~1.3.5.
The statement that $I$ extends is immediate, holding for any
isometry 
of the completion of a metric space. 
Denote the extension again by $I$.
The proof of Lemma~1.3.4  shows that the frontier
$A_{g,n}$ 
      must be
      mapped
      to itself.
      Since every point of $A_{g,n}$ is a limit of points of
$T_{g,n}$ it follows immediately that 
      the
      extension $I$ must map $A_{g,n}$ isometrically onto
itself.

      Now we wish to show the second statement; that   
     each point in $\SO(v_k)$ is mapped to  $\SO(w_k)$ for some 
$k$-simplex $w_k$.    We first show 
      inductively that  
      it is not possible for either $I$ or $I^{-1}$ 
to map a point in  $\SO(v_k)$   to a point in $\SO(w_l)$ for
 a  simplex $w_l$ with $l<k$. Since $A_{g,n}$ is mapped to
itself, the induction statement is
      true for  $k=0$. (Here $k=-1$ corresponds to $T_{g,n}$)  
Suppose the induction step for both $I$ and $I^{-1}$  is true for
all $l<k$  but a  point $x_0\in \SO(v_k)$  is mapped 
by $I$ to $y_0\in\SO(w_l)$ for $l<k$.        For some $l$-simplex
$v_l$, the point $x_0$ is in the frontier of $\SO(v_l)$.
Let $U$ be the
intersection of  a neighborhood of $x_0$ in $\ov \tg$
      with $\SO(v_l)$.  By the  induction hypothesis, the
      points in $U$
       must map into a neighborhood $V$ of $y_0$ in $\SO(w_l)$
(and not to to points in a higher dimensional space corresponding
to a lower dimensional simplex). 
      
But now $I$ is an isometry of Weil-Petersson metrics from
$U$ to $V$; consider a
      sequence in $U$ going to the frontier point $x_0$ whose
image
sequence converges
to a point in $V$.
      This again contradicts the fact that the 
      scalar
      curvature goes to $-\infty$ along the sequence
in the domain and stays bounded
on the image sequence. Thus the induction step holds for $I$. The
argument for $I^{-1}$ is identical.   
We conclude each point in $\SO(v_k)$ goes to a point in 
$\SO(w_k)$. 
Now since each $\SO(v_k)$ is a connected component in the union
(over all $k$-simplices $w_k$) of  $\SO(w_k)$,  
and $I$ is a continuous map of $A_{g,n}$, 
we  conclude that $\SO(v_k)$ must map to a single $\SO(w_k)$.

     \qed

      \section{Proof of Theorem A}

      By Corollary~1.3.6, the isometry $I$ induces an 
automorphism  $\wh I$ of the curve
      complex
      $C(F)$. By Ivanov's theorem (as extended
by Korkmaz and Luo), for $(g,n)\neq (1,2)$, the automorphism
      $\wh I$ agrees with  the automorphism of $C(F)$ induced
      by a
mapping class, say $\vp:F\to F$, which may be orientation
reversing.  As we have seen,  
the automorphism $\vp$ induces a \wpet\ isometry
$\Phi$ of $\tg$;  then  $I\circ \Phi^{-1}$ is a \wpet\ isometry
of $\tg$ extending to an isometry of $\ov T_{g,n}$ and 
preserving each frontier space $\SO(k)$.  
Thus the isometry $I\circ \Phi^{-1}$ 
induces the identity on $C(F)$. At this point then we
lose no
      generality in the argument while simplifying the notation
if
      we assume that $I=\Phi$  preserves each $\SO(v_k)$ 
(inducing the identity on $C(F)$) and we
are seeking to
      prove that $I$ is also the identity map on \tec space.
      Roughly then, our goal is to move from rough knowledge of
      $ I:\otg\to\otg$ on the frontier $A_{g,n}$ to precise
control
      on $I:\tg\to\tg$ on the interior $\tg$ of $\otg$.
 In fact we prove, by induction on the dimension 
of the frontier spaces,   that  $ I$
is the identity on each space $\SO(v_k)$ and hence on $\tg$. 

The induction hypothesis  holds for the lowest dimensional
frontier
spaces, namely the maximally pinched frontier spaces 
corresponding to the maximal 
 simplices in
$C(F)$, since the maximally pinched  frontier spaces 
are singletons and are fixed by $I$.  Our induction
hypothesis is then that for all $l>k$,  
the isometry $I$ is the identity for
all frontier Teichmuller spaces $\SO(v_l)$ of dimension $3g-4+n-
l$ and $\SO(v_k)$  is a frontier space of dimension
$3g-4+n-k$ (which is $\tg$ itself, if $k=-1$) 
Consider now $I$ restricted to $\SO(v_k)$. 
The space $\SO(v_k) = T(S)$ is the 
Teichmuller space $T(S)$ of some surface
$S$.

Let      
 $Fix=Fix(\SO(v_k))$ denote the fixed-point set of $I$ acting on
$\SO(v_k)$.   
It is a general fact that the fixed-point set of an isometry of 
a Riemannian CAT(0) space is a totally geodesic submanifold.
(To see this in our setting, note that
Wolpert's result \cite{Wol87} says that there is a geodesic 
between any two points of $\SO(v_k)$ and the negative sectional
 curvature of $\SO(v_k)$ says that this geodesic is unique. 
 Consequently, the geodesic
joining any
two points of $Fix$ is also fixed by $I$ and $Fix$ is a
convex subspace of $\SO(v_k)$. 
Next consider the action of $dI$ on the tangent space $T_zT(S)$
for $z\in Fix$.  We see that $dI$ fixes the initial tangent
vector $V_\zeta$ to any geodesic from $z$ to any other point in
$Fix$. Thus a neighborhood of $z$ in $Fix$ is given by the
exponentiated image of the kernel of $dI$. This shows that $Fix$
is a submanifold. 
Thus  
$Fix$ is a totally geodesic
submanifold of $\SO(v_k)$.)

Now suppose 
$Fix$ is a proper subset of $\SO(v_k)$. Let  $N^1(Fix) \subset
TT(S)$
be the 
unit normal bundle to $Fix$, thought of as a subbundle to 
the tangent bundle to $\SO(v_k)$.

We postpone the proof that $Fix \ne \emptyset$, while
we discuss properties of $Fix$ that follow directly
from $d_{WP}|_{\SO(v_k)}$ being $CAT(0)$. 
      
 We say that a vector $v \in
N^1_p (Fix)$
      ``exponentiates'' to  the frontier   if the 
     geodesic $\{exp_p sv|s \in [0,L]\}$
determined by $v$ joins $p$ to a frontier point
$y=exp_p Lv$. 
We claim  that in fact no vectors exponentiate to the frontier.
For
again by the induction hypothesis, if such a geodesic would
exist, 
the isometry  $I$ would fix $y$ and $p$ and
so would fix the geodesic. However, the geodesic is normal to
$Fix$,
hence not contained in $Fix$, and we have a contradiction.  

Therefore let       
    $v_0\in N^1(Fix)$ be any  vector  tangent to say
      $p_0\in Fix$, i.e. $v_0 \in T_{p_0}T(S)$.  
Let $\g_0$ the exponentiated image of
$v_0$ 
and let $p_s=\exp_{p_0}s v_0$
be a point at
      distance $s$ along $\g_0$. We claim that
\begin{claim}
$d_{WP}(p_s,Fix)=s$
\end{claim}
This claim of course follows if we know that the 
distance from $p_s$ to $Fix$ is minimized at $p_0$.  
To see this, let $q \in Fix$ be any point in $Fix$ 
other than $p_0$. The  triple $p_s, p_0, q$ of points
form a right triangle in the CAT(0) space
$T(S)$. (See \cite{BH}, Chapter II.1 for a discussion
of CAT(0) and Riemannian angles, and their
equivalence in this case.) But then the 
distance $d_{WP}(q,p_s)$ is at least as large as
the comparison Euclidean distance 
$d_{\Bbb{E}^2}(\overline{q}, \overline{p_s})$
(in the obvious notation) and
$d_{\Bbb{E}^2}(\overline{q}, \overline{p_s}) >
d_{\Bbb{E}^2}(\overline{p_0}, \overline{p_s})
= d_{WP}(p_0,p_s)$: here the inequality 
follows because $\overline{q}\overline{p_s}$
is the hypotenuse of a Euclidean right triangle,
and the equality follows by construction of the
Euclidean comparison triangle.
\qed

Note the claim says that there are points in $\SO(v_k)$
arbitrarily far
from $Fix$.  
However we will now also show that 
\begin{claim}

There is some $M$ such that {\it every}
  point in $\SO(v_k)$ is within distance $M$ of $Fix$.  
\end{claim}

This contradiction between the claims will then 
 show that $Fix=\SO(v_k)$, completing the
proof of the induction step.
 
To show that every point of $\SO(v_k)$ is within some 
distance $M$ of $Fix$, 
we first note that every point of $\SO(v_k)$ is within 
some universally bounded distance $d_0$ of {\it some} 
maximally pinched frontier point.
That statement follows from two others. The first 
(\cite{Be}, Theorem XV) is that for any
hyperbolic surface there is a maximal set of disjoint 
curves with universally bounded hyperbolic lengths.  
Such a surface is then of 
universally bounded distance from the corresponding 
maximally pinched frontier point 
by the proof in \cite{Wol75}. (Specifically,
by \cite{Str} there exists a holomorphic
quadratic differential $\Psi$ with closed 
horizontal trajectories homotopic to the elements
of a maximal set of disjoint curves such that the
corresponding cylinders have equal moduli.  Now our
upper bound for the hyperbolic lengths implies
an upper bound for the extremal lengths; it is then 
easy to check that these conditions imply a lower bound
on the moduli of those cylinders. Wolpert's proof
in \cite{Wol75}
then gives a {\it uniform} upper bound on the
Weil-Petersson length of the Teichmuller geodesic,
corresponding to $\Psi$ and tending to the maximally
pinched frontier point corresponding to the surface
pinched along each of the specified curves in the
maximal family.)

We are reduced to showing the following 
proposition.

\begin{prop}
There exists $M'$ such that  the $M'$-neighborhood of any 
maximally pinched 
frontier point $\SO(v_{3g-4+n})$ of $\SO(v_k)$ contains points
of
$Fix$. 
\end{prop}

\noindent {\bf Remark.} This claim would follow immediately,
if we knew that the interiors of geodesics between maximally
pinched frontier points  lay in a single 
component of $\overline{\SO(v_k)}$. For if so, then
consider a maximally pinched 
frontier point $\SO(v_{3g-4+n})$ of $\SO(v_k)$, and choose 
another maximally pinched 
frontier point $\SO(v'_{3g-4+n})$ of $\SO(v_k)$ so that the curves in 
$v_{3g-4+n}$ together with those in $v'_{3g-4+n}$ fill the 
surface $S$. Then a  geodesic 
between $\SO(v_{3g-4+n})$ and $\SO(v'_{3g-4+n})$ is
in $\SO(v_k)$ as well as in $Fix$, the latter because
geodesics in $CAT(0)$ spaces like $\overline{\SO(v_k)}$
are unique. As the geodesic obviously limits on 
$\SO(v_{3g-4+n})$, the Proposition follows.

\begin{proof} 
For each $x\in\SO(v_k)$ we are interested in its $I$ orbit. Let 
\begin{displaymath}
Orb(x)=\{I^j(x)\}_{j=-\infty}^\infty
\end{displaymath}

Since $I$ is an isometry
and fixes the frontier of $\SO(v_k)$  by the induction
hypothesis,
each point of  $Orb(x)$ is the same distance from the frontier
of $\SO(v_k)$  
as $x$ is. 

Fix a small $\rho>0$ and  let 
\begin{displaymath}
\SO(v_k,\rho)=\{x\in\SO(v_k):d_{WP}(x,Fr(\SO(v_k))>\rho\}
\end{displaymath}
It is easy to see $\SO(v_k,\rho)$ is open and connected. 
Let $C=C(v_{3g-4+n})$ denote the curves in 
$v_{3g-4+n}\setminus v_k$. 
Now let 
\begin{displaymath}
 \Omega(C,\rho)=\{x\in\SO(v_k,\rho): 
\sup_{y\in Orb(x)} l_C(y)<\infty\}
\end{displaymath}
That is, $\Omega(C,\rho)$ consists of those points such that the
lengths of the curves in $C$ are bounded on  the entire orbit. 

We begin by claiming that 

\begin{claim} $\Omega(C,\rho)=\SO(v_k,\rho)$.
\end{claim}

We prove the claim by showing that $\Omega(C,\rho)$
is non-empty, open and closed in the connected set
$\SO(v_k,\rho)$.

We first show that $\Omega(C,\rho)$ is nonempty. 
Using the isometric action of the mapping class group there
exists $\epsilon>0$ and $\rho>0$ and $x_0$ so that
$d_{WP}(x_0,\SO(v_{3g-4+n}))=\epsilon$
 and
$d_{WP}(x_0,Fr(\SO(v_k))=2\rho$. 
We may choose $\epsilon$ (and $\rho$)
small enough so that the hypothesis of
Lemma~1.3.2  holds.  Then since the orbit of $I$ remains within
$\epsilon$ of $\SO(v_{3g-4+n})$, Lemma~1.3.2 says that $x_0\in
\Omega(C,\rho)$.

We now show that $\Omega=\Omega(C,\rho)$ is open. 
Let $x\in\Omega$  and let $L= L(x) =\sup_{y\in Orb(x)} l_C(y)$.
  By 
Lemma~1.3.3, there exists $\delta=\delta(\rho,L)$  such that for
all  $y\in
Orb(x)$ and all $\zeta$ in the $\delta$ ball about $y$, 
\begin{displaymath}
l_C(\zeta)\leq 2L
\end{displaymath}
   Since
$I$ is an isometry, the $I$ orbit of a $\delta$ ball about $x$ is
the union of the $\delta$ balls about the points in $Orb(x)$. 
Since $\SO(v_k,\rho)$ is open, if we take  $\delta$ small enough
we can insure the $\delta$ ball about $x$ remains in
$\SO(v_k,\rho)$ and the inequality above says that it is then
contained in $\Omega$, proving $\Omega$ is open. 

Finally, we show $\Omega$ is closed. 
Let $z_0$ be a limit of points $z_i\in\Omega$ and assume
$z_0\notin\Omega$. Then there exists a sequence
$\{y_j=I^j(z_0)\} \subset
Orb(z_0)$ such that $l_C(y_j)\to\infty$.  Since
$y_j\in\SO(v_k,\rho)$  there is a compact set
$K\subset\SO(v_k,\rho)$ and a sequence $\{f_j\} \subset Mod(g,n)$ such
that $\zeta_j=f_j(y_j)\in K$. Let $C_j=f_j(C)$ be the collection
of
image curves. We have
\begin{displaymath}
l_{C_j}(\zeta_j) = l_{f_j(C)}(f_j(y_j)) = l_C(y_j) \to\infty
\end{displaymath}

Since $I$ is an isometry and $f_j$ induces an isometry, we have
that
 for $i$ large enough, the points $\zeta_{j,i}=f_j(I^j(z_i))$ lie
in $K$.  Fix such an index $i$ and consider the sequence
$\{\zeta_{j,i}\}$.  Now we have  
\begin{displaymath}
\sup_j l_{C_j}(\zeta_{j,i})<\infty
\end{displaymath}
However,  given a compact set $K$, there exists a constant $L_0$
(depending only on $K$)
such that for any two points $x,y\in K$ and {\em any} curve
$\beta$ we have $\frac{l_\beta(x)}{l_\beta(y)}\leq L_0$.  This
contradicts the previous assertions 
that the ratio of $l_{C_j}(\zeta_j)$ and $l_{C_j}(\zeta_{j,i})$ 
have no bound; thus $z_0$ must in
fact lie in $\Omega$, and $\Omega$ is closed in $\SO(v_k,\rho)$.
This concludes the proof of the claim.
 
We now conclude the proof of the Proposition.
Now let   
$\SO(w_{3g-4+n})$ be another maximally pinched frontier
point
of $\SO(v_k)$ such that the curves in $C=v_{3g-4+n}\setminus v_k$
together with those in 
$C'=w_{3g-4+n}\setminus v_k$   {\em fill}
 the surface $S$,  which  means
that if we remove the curves in $C\cup C'$ from $S$, the result is a 
union of simply
connected domains. 

We  form the corresponding  $\Omega(C',\rho)$, and again 
the claim above shows that this
coincides with $\SO(v_k,\rho)$.  Thus every point in $\SO(v_k,\rho)$ lies in
  $\Omega(C,\rho)\cap\Omega(C',\rho)$.

Let $x_0$ be a point at distance $2\rho$ from the frontier and $\epsilon$ from 
 $\SO(v_{3g-4+n})$. Applying Lemma~1.3.3 again, 
there is a ball $B$ of radius $\delta$ about
$x_0$ such that if we set  
\begin{displaymath}
B_0= \cup_{j=-\infty}^{\infty}I^j(B)
\end{displaymath} 
then there exists $M_1$  such that 
 $l_{C\cup C'}(x)<M_1$ for all $x\in B_0$. 
We may take $\delta<\rho$. 
Clearly $B_0$ is $I$ invariant. 

Now we  adapt an argument of Wolpert's
\cite{Wol87}
 in his proof of the Nielsen Realization Problem that every
finite subgroup of the mapping class group has a fixed point
to find a fixed point in the current situation.  
 
Consider the  subset 
\begin{displaymath}
W=W(C\cup C')=\{x\in \SO(v_k):l_{C\cup C'}(x)\leq M_1\}
\end{displaymath} 
By \cite{Wol87},  
since the set of curves $C\cup C'$ fills $S$,
the set $W$ is a compact 
subset of $\SO(v_k)$; also the set $W$  is
a 
cell, its  boundary is $C^1$, and it contains $B_0$ in its
interior. 

Now define a function $D$ on $\SO(v_k)$ by 
\begin{displaymath}
D(x)=\frac{1}{\mu(B_0)}\int_{B_0} d_{WP}(x,y) d\mu(y), 
\end{displaymath}
where $\mu$ is  Weil-Petersson volume element. 
As $B_0$ is a subset of the compact set $W$,
the set $B_0$ has  finite total measure.  Further,  as
$d_{WP}(x,\cdot)$ is bounded on $W \supset B_0$, we see that
$D(x)$ is a finite integral for each $x$.

Since $B_0$ is $I$-invariant, so is $D$.  Since
$l_C(x)>l_C(y)$ for $x$ in the boundary of $W$ and $y\in B_0$,  
the strict convexity of $l_C(x)$ along Weil-Petersson geodesics
says that the vector 
 $\grad_{WP}d_{WP}(\cdot,y)$ points out at each
point on the boundary of $W$.  Thus as an average of such
vectors, 
$\grad_{WP}D(\cdot)$ points
out at each such boundary point. By the Poincare-Hopf index
theorem, we see that there is some $y_0\in W$
for which $\grad_{WP}D(y_0)=0$.  
  By the negative curvature and the geodesic
 convexity of the metric, the distance from a point to a geodesic
is a strictly convex function of the parameter along the
geodesic.  
Consequently as an average of 
such functions, $D$ is also strictly convex. This implies that
$y_0$ is the unique
minimum for $D_0$. It follows that 
$I(y_0)=y_0$.

Now we wish to estimate $d_{WP}(y_0,\SO(v_{3g-4+n}))$.  Each point of $B_0$ 
is distance at most $\epsilon+\delta\leq \epsilon+\rho$
 from $\SO(v_{3g-4+n})$. Let 
$y_1\in\SO(v_k)$ be any point within distance $\rho$ of $\SO(v_{3g-4+n})$ and
 thus  $y_1$ is at most $\epsilon+2\rho$ from any point of $B_0$.
  Since $D$ is the average of such distances, $D(y_1)\leq \epsilon+2\rho$.
Since $y_0$ is the  point that minimizes $D$, \begin{displaymath}
D(y_0)\leq D(y_1)\leq \epsilon+2\rho
\end{displaymath} 
Since $D(y_0)$ is the  average of distances frm $y_0$ to 
points in $B_0$, 
\begin{displaymath}
\min_{y\in  B_0} d_{WP}(y_0,y)\leq \epsilon+2\rho.
\end{displaymath}
Since we have the bound $\epsilon+\rho$ on the distance from any point
 of $B_0$ to $\SO(v_{3g-4+n})$,
 \begin{displaymath}
d_{WP}(y_0, \SO(v_{3g-4+n}))\leq 2\epsilon+3\rho
\end{displaymath}
We have thus 
found a point of $Fix$ within uniform distance $M'=2\epsilon+3\rho$
 of every
maximally pinched frontier  of $\SO(v_k)$.
\end{proof} 
The proof of the main Theorem A is now complete.

\noindent
Howard Masur:\\
Dept. of Mathematics, University of Illinois, Chicago\\
851 S. Morgan\\
Chicago, Il 60607\\
E-mail: masur@math.uic.edu
\medskip

\noindent
Michael Wolf:\\
Dept. of Mathematics, Rice University\\
Houston, TX 77251\\
E-mail: mwolf@math.rice.edu


\begin{thebibliography}{MW}

\bibitem {Ab71} W.Abikoff, {Augmented Teichmuller spaces},
{\em Bull. Amer. Math. Soc.} {\bf 82} (1971) 333-334

\bibitem {Ab77} W.Abikoff, {Degenerating families of Riemann
surfaces}
{\em Annals of Math.} {\bf 105} 1977 29-44

\bibitem {Ahl} L. Ahlfors, {Some Remarks on Teichm\"uller's
Space of Riemann Surfaces} {\em Ann. Math.} {\bf 74}(1961)
171-191

\bibitem{Ba} W. Ballman, {\em Lectures on spaces of nonpositive
curvature}
Birkhauser Verlag , Basel (1995)


\bibitem {Be} L.Bers, {Spaces of degenerating Riemann surfaces in
Discontinous groups and Riemann surfaces}, {\em Annals of Math
Studies} {\bf 79} Princeton
University, Princeton New Jersey, (1974)

\bibitem{BH} M.Bridson, A.Haefliger, {Metric Spaces of Non-Positive 
Curvature}, 
Springer Berlin (1999)


\bibitem{Chu75} T. Chu, {The Weil-Petersson metric in moduli
space}
{\em Chinese J. Math.} {\bf 4} (1976) 29-51

\bibitem{EK} C.Earle, I.Kra,
{On isometries  between Teichmuller spaces}, {\em Duke J. Math}
{\bf 41} (1974)
583-591


\bibitem {Iv97} N. Ivanov, {Automorphism of complexes of curves
and of Teichmuller spaces},
{\em Internat. Math. Res. Notices} {\bf 14} (1997) 651-666

\bibitem {Kor97} M. Korkmaz, {Automorphisms of Complexes of
Curves in 
Punctured Spheres and on Punctured Tori}, {\em Top. and its
Appl.},
to appear

\bibitem{Luo00} F. Luo, {Automorphisms of the complex of curves}
{\em Topology}
{\bf 39} (2000) 283-298
 
\bibitem{Ma76} H.Masur, {The extension of the Weil-Petersson
metric to the boundary of
Teichmuller space}, {\em Duke Math J}. {\bf 43} (1976) 623-635

\bibitem{Mum} D. Mumford {A remark on Mahler's compactness theorem},
{\em Proc. Amer. Math. Soc.} {\bf 28}(1971) 289-294


\bibitem {Nag} S.Nag, {Complex Analytic theory of Teichmuller
space}, Wiley, New York, (1988)

\bibitem {Ro} H. Royden, {Automorphisms and isometries of
Teichmuller spaces. Advances in
the
theory of Riemann surfaces}, {\em Annals of Math. Studies}, {\bf
66} (1970) 369-383

\bibitem {Roy74} H. Royden, Oral communication.

\bibitem {SY78} R.Schoen, S.-T. Yau, {Compact group actions and
the topology of manifolds with
nonpositive curvature} {\em topology} {\bf 18} (1979) 361-380

\bibitem{Str} K. Strebel {Quadratic Differentials},
Springer, Berlin, (1984)

\bibitem {Tra92} S. Trapani, {On the determinant of the bundle of
meromorphic 
quadratic differentials on the Deligne-Mumford compactification
of the moduli space of Riemann
surfaces} {\em Math. Ann.} {\bf 293} (1992) 681-705

\bibitem{Tro86} A. Tromba {On a natural algebraic affine connection
on the space of almost complex structures and the curvature 
of Teichmuller space with respect to its Weil-Petersson metric}
{\em Manuscripta Math.} {\bf56}(1986) 475-497

\bibitem {Tro92} A.Tromba,
{Teichmuller theory in Riemannian geometry} {\em Birkhauser}
Basel
(1992)

\bibitem{Wol75} S. Wolpert,
{Noncompletenes of the Weil-Petersson metric for Teichmuller
space}. {\em Pacific J.Math}. {\bf
61} (1975)
573-577

\bibitem{Wol76} S. Wolpert, {The finite Weil-Petersson diameter of 
Riemann space} {\em Pacific J. math.} {\bf 70} (1977) 281-288

\bibitem{Wol86} S. Wolpert, {Chern forms and the Riemann 
tensor for the moduli space of curves}. {\em Invent. Math.}
{\bf85}(1986) 119-145

\bibitem {Wol87} S.Wolpert, {Geodesic length functions and the
Nielsen problem},
{\em J.Diff. Geom.} {\bf 25} (1987) 275-295


\bibitem {Y} S.Yamada, {Weil-Petersson Completion of 
Teichmuller spaces and mapping class
group actions}, manuscript

\end{thebibliography}
      \end{document}